\documentclass[a4paper,11pt,reqno, english]{amsart}

\usepackage{anysize,color}
\usepackage[utf8]{inputenc}
\usepackage[dvipsnames]{xcolor}
\usepackage[colorlinks=true,urlcolor=blue,linkcolor={YellowOrange},citecolor={YellowOrange}]{hyperref}  
\usepackage{float}
\usepackage{amsfonts,amssymb,amsmath}
\usepackage{amsthm}    
\usepackage{url}
\usepackage{paralist}
\usepackage{tikz}
\usepackage{subcaption}
    \usetikzlibrary{decorations.markings,arrows,shapes.callouts}
\usepackage{rotating} 
\usepackage[mathscr]{euscript}
\usepackage{verbatim}
\usepackage{tikz-cd}
\usepackage{mathtools}
\usepackage{dsfont}
\usepackage{amsrefs}
\usepackage[capitalise]{cleveref}
\usepackage{thm-restate}

\usepackage{charter}

\newtheorem{theorem}{Theorem}[section]

\newtheorem{lemma}[theorem]{Lemma} 
\newtheorem{conjecture}[theorem]{Conjecture}

\theoremstyle{definition}
\newtheorem{definition}[theorem]{Definition}
\theoremstyle{remark}

\newcommand{\rr}{\mathds{R}}
\newcommand{\R}{\mathds{R}}

\newcommand{\PP}{\mathbb{P}}

\newcommand{\Z}{\mathbb{Z}}

\newcommand{\FF}{\mathbb{F}}

\newcommand{\FHindex}[3]{\mathrm{ind}^{#1}_{#2}(#3)}

\newcommand{\grassR}[2]{G_{#1}(\R^{#2})}

\DeclareMathOperator{\id}{\mathrm{id}}
\DeclareMathOperator{\pt}{\mathrm{pt}}

\tikzcdset{row sep/my_size/.initial=-2mm}
\tikzcdset{row sep/my_size_small/.initial=-3mm}
\newcommand{\xDownarrow}[1]{%
  {\left\Downarrow\vbox  to #1{}\right.\kern-\nulldelimiterspace}
}

\begin{document}

\title{Bisections of mass assignments by parallel hyperplanes}

\author[Sadovek]{Nikola Sadovek} 
\address{Institute of Mathematics, Freie Universit\"at Berlin, Germany}
\email{nikolasdvk@gmail.com,~nikola.sadovek@fu-berlin.de}

\author[Sober\'on]{Pablo Sober\'on} 
\address{Baruch College and The Graduate Center, City University of New York, New York, USA.}
\email{psoberon@gc.cuny.edu}

\thanks{The research of N. Sadovek is funded by the Deutsche Forschungsgemeinschaft (DFG, German Research Foundation) under Germany's Excellence Strategy – The Berlin Mathematics Research Center MATH+ (EXC-2046/1, project ID 390685689, BMS Stipend). The research of P. Sober\'on was supported by NSF CAREER award no. 2237324 and a PSC-CUNY Trad B award.}

\begin{abstract}
    In this paper, we prove a result on the bisection of mass assignments by parallel hyperplanes on Euclidean vector bundles. Our methods consist of the development of a novel lifting method to define the configuration space--test map scheme, which transforms the problem to a Borsuk--Ulam-type question on equivariant fiber bundles, along with a new computation of the parametrized Fadell--Husseini index. As the primary application, we show that any $d+k+m-1$ mass assignments to linear $d$-spaces in $\mathbb{R}^{d+m}$ can be bisected by $k $ parallel hyperplanes in at least one $d$-space, provided that the Stirling number of the second kind $S(d+k+m-1, k)$ is odd. This generalizes all known cases of a conjecture by Sober\'on and Takahashi, which asserts that any $d+k-1$ measures in $\mathbb{R}^d$ can be bisected by $k$ parallel hyperplanes.
\end{abstract}

\maketitle



\section{Introduction}

\subsection*{Motivation} Problems concerning the existence of equitable distributions of measures in $\rr^d$ using hyperplane arrangements have motivated significant developments in discrete geometry and topological combinatorics \cites{RoldanPensado2022, Zivaljevic2017}.  These problems ask whether any given collection measures can be fairly partitioned by certain hyperplane arrangements.  The ham sandwich theorem of Steinhaus \cite{Steinhaus1938}, one of the most emblematic results in mass partitions, fits this description. It states that \emph{any $d$ finite absolutely continuous measures on $\rr^d$ there exists a hyperplane that simultaneously halves all of them.} Here, by an \textit{absolutely continuous} measure on $\rr^d$ we mean absolutely continuous with respect to the Lebesgue measure. 

As we impose additional constraints in the family of hyperplanes -- such as parallelism, orthogonality, or fixed intersection patterns -- the question of which families of measures can be fairly partitioned by them quickly brings the problem to the forefront of current research (see, e.g. \cites{Blagojevic:2018jc, Hubard2020, Aronov2024}). In particular, Sober\'on and Takahashi conjectured the following.

\begin{conjecture}[\cite{Soberon2023}]
    \label{conj:soberon-takahashi}
    Let $d$ and $k$ be positive integers.  Then, any $d+k-1$ finite absolutely continuous measures on $\R^d$ can be simultaneously bisected by the chessboard coloring induced by $k$ or fewer parallel hyperplanes.
\end{conjecture}

Given a set of parallel hyperplanes, the \textit{chessboard coloring} they induce is the result of coloring $\rr^d$ with two colors, white and black, alternating color as we go through a hyperplane.  We consider the white and black sets as closed: their common boundary is the union of the family of hyperplanes, which has Lebesgue measure zero, so it causes no problem.

The case $k=1$ of Conjecture \ref{conj:soberon-takahashi} is the ham sandwich theorem, and the case $d=1$ is the necklace splitting theorem for two thieves \cites{Hobby:1965bh, Goldberg:1985jr}. Furthermore, the case $k=2$ of Conjecture \ref{conj:soberon-takahashi} was originally proved by Sober\'on and Takahashi \cite{Soberon2023} and later obtained by Blagojevi\'c and Crabb \cite{Blagojevic2025} via different methods. Recently, Hubard and Sober\'on extended these results as follows:

\begin{theorem}[\cite{HubardSoberon24}]
\label{thm:hubard-soberon} 
Let $k,d$ be positive integers such that the Stirling number of the second kind $S(d+k-1,k)$ is odd.  Then, for any $d+k-1$ absolutely continuous finite measures on $\rr^d$ there exists a family of $k$ or fewer parallel hyperplane such that the chessboard coloring they induce bisects each measure.
\end{theorem}


 Another family of mass partition problems that has been the focus of recent research is the partition of mass assignments.  Given integers $m$ and $d$, a \emph{mass assignment} is a function that continuously assigns a smooth finite measure to each $d$-dimensional subspace of $\rr^{d+m}$. In part, the following observation motivates the consideration of mass assignment partitions. Namely, the ham sandwich theorem is optimal with respect to the number of measures in the sense that for $d+1$ finite absolutely continuous measures in $\R^d$ there does not have to exist a bisecting hyperplane, as is evident by an example of measures where each is placed on ball of small radius around one of $d+1$ affinely independent points in $\R^d$. However, a natural question arises: \emph{given $d+1$ mass assignments to linear $d$-subspaces of $\R^{d+1}$, can one find a linear $d$-space $V$ such that the measures assigned to it can be simultaneously bisected by a hyperplane in $V$?}
Recently, it has been shown that for many mass partition problems, the additional freedom of being able to choose a $d$-dimensional subspace indeed allows us to bisect more measures than in the classic result.  The primary example is the ham sandwich for mass assignments by Schnider:

\begin{theorem}[\cite{Schnider:2020kk}]
    \label{thm:schnider}
    For any $d+m$ mass assignments to the $d$-dimensional linear subspaces of $\rr^{d+m}$, there exists such a $d$-dimensional subspace $V$ and a hyperplane in $V$ that halves each of the $d+m$ measures assigned to it.
\end{theorem}

In contrast, a direct application of the ham sandwich theorem would only allow us to bisect $d$ measures in $V$.  
Mass assignment problems for partitions with hyperplane arrangements lead to new questions in equivariant topology, such as results on non-existence of equivariant maps between equivariant fiber bundles \cite{Blagojevic2025}. Thus, novel methods are required to study corresponding configurations spaces and their symmetries. Recent examples of results on partitions of mass assignments include the fairy bread sandwich theorem \cites{AxelrodFreed2024}, a mass assignment versions of the Gr\"unbaum-Hadwiger-Ramos problem \cites{Blagojevic2025, Blagojevic2023}, and a mass assignment version of the generalized Nandakumar--Ramana-Rao problem \cite{levinson2023convex}.

As with the ham sandwich theorem, \cref{conj:soberon-takahashi} could not be improved in the number of measures: $d+k$ measures, each placed on a ball of small radius around one of $d+k$ points in general position, cannot be bisected by $k$ parallel hyperplanes. This raises the question of the existence of bisections of mass assignments by chessboard colorings induced by parallel hyperplanes.

\subsection*{Our results} In this manuscript, we consider a mass assignment extension of Conjecture \ref{conj:soberon-takahashi}. The first main result is the following generalization of \cref{thm:hubard-soberon}.

\begin{restatable}{theorem}{thmMain}
\label{thm:assignment}
    Let $d \ge 1$, $k \ge 1$ and $m \ge 0$ be integers such that the Stirling number of second kind $S(d+m+k-1,k)$ is odd.  For any $d+m+k-1$ mass assignments to $d$-dimensional linear subspaces of $\R^{d+m}$ there exists a $d$-dimensional linear subspace $V$ of $\rr^{d+m}$ and $k$ or fewer hyperplanes in $V$ whose induced chessboard coloring bisects each of the measures assigned to $V$.
\end{restatable}

Indeed, the case $m=0$ of Theorem \ref{thm:assignment} recovers \cref{thm:hubard-soberon} with exactly the same condition on the parameters, while the case $k=1$ recovers Theorem \ref{thm:schnider}, as presented in the diagram:
\vspace{-4mm}
\begin{figure*}[h]
    \centering
    \begin{equation*}
    \begin{tikzcd}[row sep = -4.5mm, column sep = 1.2mm]
        \text{Theorem \ref{thm:assignment}} & & \xRightarrow{(m=0)} & & \text{Theorem \ref{thm:hubard-soberon}}\\
         & & & & \text{(Hubard \& Sober\'on)} \\
        \phantom{\xDownarrow{3mm}} & & & & {}\\
        \hspace{7mm}\xDownarrow{4mm} \overset{(k=1)}{} & &  & &  \xDownarrow{4mm} \overset{(k=1)}{} \\
        \phantom{\xDownarrow{3mm}} & & &  & \\
        \text{Theorem \ref{thm:schnider}} & & \xRightarrow{(m=0)} & & \text{Ham sandwich theorem} \\
        \text{(Schnider)} & & & & \text{(Steinhaus)}\\
    \end{tikzcd}
\end{equation*}
    \label{fig:implications}
\end{figure*}

\vspace{1mm}

\noindent Moreover, we note that the case $k=2$ of Theorem \ref{thm:assignment} was obtained independently by Lessure and Sober\'on \cite{LessureSoberon25} using different methods. On the other hand, the case $d=1$ of Theorem \ref{thm:assignment} is a curious extension of the necklace splitting theorem.  Even though we do not have a natural interpretation in terms of thieves and necklaces, \cref{thm:assignment} provides insight into the kind of degrees of freedom that are required on a family of $m+k$ measures on $\rr^1$ in order to bisect all of them using at most $k$ cuts.

Our second main result is more general Theorem \ref{thm:mass-assign-general} regarding bisections of mass assignments on Euclidean vector bundles. There, the parity condition of the Stirling number of second kind is replaced by a more general algebraic condition of an ideal non-containment, expressed in terms of Stiefel-Whitney classes of the bundle. Theorem \ref{thm:assignment} is then a special case of Theorem \ref{thm:mass-assign-general} applied to the canonical bundle over the Grassmann manifold of $d$-planes in $\R^{d+m}$.

The methods that Hubard and Sober\'on used to prove \cref{thm:hubard-soberon} are purely geometric, based on a technique by Schnider \cite{Schnider2021}.  In contrast, our proof of \cref{thm:assignment} is topological, based on a new lifting method used to define the configuration space -- test map scheme and novel computations of the Fadell--Husseini index of equivariant fiber bundles. The fact that the same parity condition arises with two completely different methods may indicate that it could be necessary.

\subsection*{Organization of the paper} In Section \ref{sec:prelim} we introduce concepts and recall classical results used throughout the paper. In Section \ref{sec: cs-tm scheme I} we develop a new lifting method and the configuration space -- test map scheme for a more general problem of chessboard coloring bisections of mass assignments on Euclidean vector bundles. In Section \ref{sec:equiv-bundles} we compute Fadell-Husseini indices of relevant equivariant bundles and prove a more general Theorem \ref{thm:mass-assign-general} regarding the existence of such mass assignment bisections. Finally, in Section \ref{sec:proof-main}, we prove Theorem \ref{thm:assignment}. 

\section{Notation and preliminaries}
\label{sec:prelim}

\subsection*{Fadell-Husseini index}
For a finite group $G$, let us denote by $EG$ a contractible free $G$-space and by $BG =EG/G$ the classifying space. For a $G$-space $X$, let $X \to X\times_G EG \to BG$
be the associated Borel fiber bundle, where $X\times_G EG$ is the quotient of $X\times EG$ by the diagonal $G$-action. Given a ring $R$, we denote by $H_G^*(X;R) = H^*(X\times_G EG;R)$ the $G$-\emph{equivariant cohomology} of $X$. We note that this construction is functorial, in the sense that a $G$-equivariant map $f\colon X \to Y$ induces a map $f \times_G \id \colon X \times_G EG \to Y \times_G EG$, and hence a morphism
\[
    f^* \colon H^*_G(Y;R) \longrightarrow H^*_G(X;R).
\]
For a fiber bundle $F \to E \xrightarrow{\pi} B$ such that $E$ and $B$ are $G$-spaces and the projection map $\pi$ is $G$-equivariant, the ideal
\begin{align*}
    \FHindex{R}{G}{\pi} =\ker \left(\pi^* \colon~ H^*_G(B;R) \longrightarrow H^*_G(E;R) \right)  \subseteq H^*_G(B;R)
\end{align*}
is the \emph{Fadell-Husseini index} of $\pi$ (see \cite{Fadell:1988tm}). If $F' \to E' \xrightarrow{\pi'} B$ is another such fiber bundle and
\begin{equation*}
    \begin{tikzcd}
        E' \arrow[r,"f"] \arrow[d, "\pi'"] & E \arrow[d, "\pi"]\\
        B \arrow[r, equal] & B
    \end{tikzcd}
\end{equation*}
is a $G$-equivariant bundle map, $\emph{monotonicity}$ of Fadell-Husseini index \cite{Fadell:1988tm} states that
\begin{equation*}
    \FHindex{R}{G}{\pi} \subseteq \FHindex{R}{G}{\pi'} \subseteq H^*_G(B).
\end{equation*}

\subsection*{Representations}
Throughout the paper, we will consider the group $G = \Z_2^2$. For $i=1,2$, let $\R_i$ denote the real 1-dimensional $\Z_2^2$-representation where the $i$-th $\Z_2$-factor acts antipodally, while the other factor acts trivially. Moreover, let us denote by $\R_{1,2}=\R_1 \otimes_{\R} \R_2$. Next, we define a particular representation that will be used in the following sections.

\begin{definition} \label{def:representation-R,a,b}
    For integers $a, b\ge 0$, let $R^{a,b} \coloneqq \R_{1,2}^{\oplus a} \oplus \R_{2}^{\oplus b}$ denote a real $\Z_2^2$-representation of dimension $a+b$.
\end{definition}

\subsection*{Grassmannians}
For integers $d \ge 1$ and $m \ge 0$, let $\grassR{d}{d+m}$ denote the Grassmann manifold of linear $d$-subspaces in $\R^{d+m}$.
Its cohomology ring \cite{borel1953cohomologieMod2} is given by
\begin{equation} \label{eq:coh-grass}
	H^*(\grassR{d}{d+m};\FF_2) \cong \FF_2[w_1, \dots, w_d,\overline{w}_1, \dots, \overline{w}_m]/I,
\end{equation}
where $|w_i| = i$ and $|\overline{w}_j|=j$, and the ideal $I$ is
\begin{equation} \label{eq:cohom-grass}
	I = \left(w_i + w_{i-1}\overline{w}_1 + \dots + \overline{w}_i \colon~ i=1,  \dots, d+m \right).
\end{equation}
Here, it is assumed that $w_{i} =  0  $, for $i >d$, and $\overline{w}_{j} = 0$, for $j > m$. We will denote by $\gamma_d$ the canonical rank $d$ vector bundle over $\grassR{d}{d+m}$. We recall  that $w_0, w_1, \dots, w_d \in H^*(\grassR{d}{d+m},\Z_2)$ are the Stiefel-Whitney classes of $\gamma_d$ (see \cites{MilnorStasheff74, hatcher2003vector}).

\subsection*{Mass assignments}
For a Euclidean vector space $V$, let $\mathcal{M}(V)$ be the space of finite measures on $V$ which are absolutely continuous with respect to the Lebesgue measure in $V$, equipped with the weak topology.
For a Euclidean vector bundle $\xi \colon \R^d \to E \xrightarrow{\pi} B$, let
\begin{equation*}
    \mathcal{M}(\xi) \colon \mathcal{M}(\R^d) \longrightarrow \{(b, \nu)\colon \in b \in B,~ \nu \in \mathcal{M}(\pi^{-1}(b))\} \longrightarrow B
\end{equation*}
be the induced measure-bundle with the projection map $(b, \nu) \mapsto b$. The topology of its total space is defined using the local triviality of $\pi$ (see \cites{Blagojevic2023, Blagojevic2025}). A section of $\mathcal{M}(\xi)$ is called a \emph{mass assignment} on $\xi$. Given one such mass assignment $\mu$, we will denote by $\mu[b] \coloneqq \mu(b) \in \mathcal{M}(\pi^{-1}(b))$ the measure that $\mu$ induces on $\pi^{-1}(b)$. 

\begin{definition} \label{def:mass-assign}
    A \emph{mass assignment to the $d$-dimensional subspaces of $\R^{d+m}$} is a mass assignment on the canonical bundle $\gamma_d$ over $\grassR{d}{d+m}$.
\end{definition}

\section{Configuration space -- test map scheme}
\label{sec: cs-tm scheme I}

In this section, we use a novel lifting method to develop a configuration space -- test map scheme for Theorem \ref{thm:assignment}. For other methods applied to math assignments, see \cite{Blagojevic2025}. For the lemma below we use notation from Section \ref{sec:prelim}.

For a Euclidean vector bundle $\xi \colon~ \R^d \to E \to B$ over a paracompact base space with a fiberwise action of $\Z_2^2$ given by the antipodal action of the first factor $\Z_2$ we have the induced fiber bundle
\begin{equation*}
    S(\R_1^d) \times S(R^{a,b}) \longrightarrow S(E) \times S(R^{a,b}) \longrightarrow B
\end{equation*}
with a fiberwise action of $\Z_2^2$. Since the bundle is Euclidean, its fiberwise scalar product \cite{hatcher2003vector}*{Prop.~1.2} defines the sphere bundle.

\begin{lemma} \label{lem:cs-tm-parallel}
    Let $d \ge 2$, $k \ge 1$ and $n \ge 1$ be integers. Suppose that $\xi \colon~ \R^d \to E \to B$ is a Euclidean vector bundle over a paracompact base space with a fiberwise action of $\Z_2^2$ given by the antipodal action of the first factor $\Z_2$ and let $\mu_1, \dots, \mu_{n}$ be mass assignments on $\xi$. Assume that there does not exists a $\Z_2^2$-equivariant map of bundles
    \begin{equation} \label{eq:bundle-equiv}
        \begin{tikzcd}
            S(E) \times S(R^{\lceil k/2 \rceil,\lfloor k/2 \rfloor}) \arrow[r] \arrow[d] &  S(\R_2^{n-1}) \times B \arrow[d]\\
            B \arrow[r, equal] & B
        \end{tikzcd}
    \end{equation}
    over the identity map on $B$, where $R^{a,b}$ is the $\Z_2^2$-representation from Definition $\ref{def:representation-R,a,b}$. 
    Then, there exists $b \in B$ such that the masses $\mu_1[b], \dots, \mu_n[b]$ in the fiber $F_b \cong \R^d$ of $\xi$ over $b$ can be simultaneously bisected by a chessboard coloring induced by $k$ or fewer parallel hyperplanes in $F_b$.
\end{lemma}

The rest of this section is dedicated to proving Lemma \ref{lem:cs-tm-parallel}. To this end, let $b \in B$ be a point and $\mu_1[b], \dots, \mu_{n}[b]$ be absolutely continuous measures in the fiber $F_b \cong \R^d$. We will construct a $\Z_2^2$-equivariant \textit{test map}
\begin{equation} \label{eq:test-map}
    f_b \colon~ S(F_b) \times S\left((\R_{1,2} \oplus \R_2)^{\lfloor k/2 \rfloor} \oplus \R_{1,2}^{k~(\text{mod}~2)}\right) \longrightarrow \R_2^{d+k-2},
\end{equation}
such that if $f_b(v,n) = 0$, then there exist $k$ or fewer parallel hyperplanes in $F_b$ with the bisecting property we seek. We note that
\[
    (\R_{1,2} \oplus \R_2)^{\lfloor k/2 \rfloor} \oplus \R_{1,2}^{k~(\text{mod}~2)} \cong \R_{1,2}^{\lceil k/2 \rceil} \oplus \R_2^{\lfloor k/2 \rfloor} = R^{\lceil k/2 \rceil,\lfloor k/2 \rfloor}
\]
as $\Z_2^2$-representation, but we will use former form to define $f_b$. Finally, our construction will be continuous with respect to $b \in B$ and will induce a $\Z_2^2$-equivariant bundle map 
\begin{equation} \label{eq:bundle-test-map}
        \begin{tikzcd}
            S(E) \times S(R^{\lceil k/2 \rceil,\lfloor k/2 \rfloor}) \arrow[r, "f"] \arrow[d] &  \R_2^{n-1} \times B \arrow[d]\\
            B \arrow[r, equal] & B
        \end{tikzcd}
\end{equation}
over the identity on $B$. The proof of the lemma is then finalized once we show that the image of $f$ intersects the zero section in the codomain bundle $\R_2^{n-1} \times B \to B$. Indeed, if that was not the case, there would exist an equivariant bundle map \eqref{eq:bundle-equiv} from Lemma \ref{lem:cs-tm-parallel} obtained by postcomposing $f$ with the equivariant retraction 
\[
    \left(\R_2^{n-1} \setminus \{0\}\right) \times B \longrightarrow S(\R_2^{n-1}) \times B,~ (v,b) \longmapsto (v/\|v\|, b).
\]
However, this is a contradiction with the assumption from the lemma. 

In the reminder of the section we will define the $\Z_2^2$-equivariant test map \eqref{eq:test-map}, for each $b \in B$, which then induces the bundle map \eqref{eq:bundle-test-map}; finally, we will show that zeros of $f_b$ yield bisections of masses $\mu_1[b], \dots, \mu_n[b]$ in $F_b$ by parallel hyperplanes.

\subsection*{Definition of the test map}

\begin{figure}
\centering
\begin{subfigure}{.4\textwidth}
  \centering
  \includegraphics[width=\linewidth]{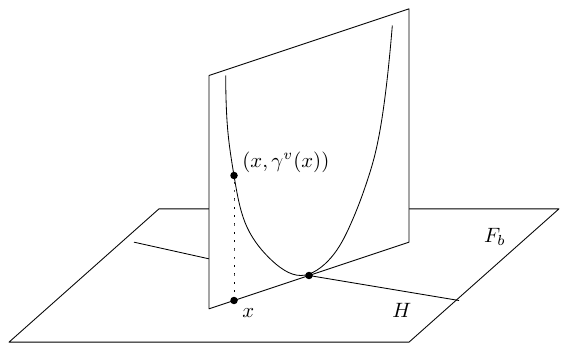}
\end{subfigure}%
\begin{subfigure}{.4\textwidth}
  \centering
  \includegraphics[width=\linewidth]{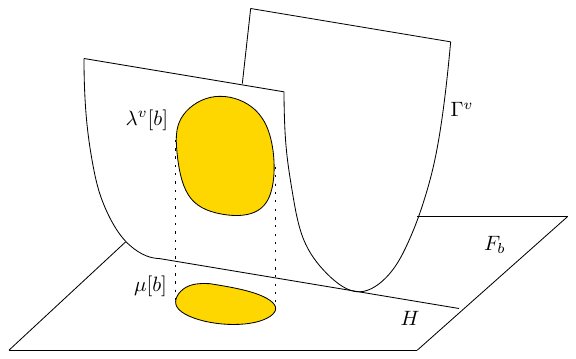}
\end{subfigure}
\caption{The map $\gamma^v$ (left) and the lift $\lambda^v[b]$ of a measure $\mu[b]$ (right) to the graph $\Gamma^v$ of $\gamma^v$.}
\label{fig:lifts}
\end{figure}

Let $b \in B$ be fixed and let $(v,n) \in S(F_b) \times S(\R^k)$. Following \cite{Soberon2023}*{Lem.~3.1}, we define an affine hyperplane $H \subseteq F_b$ parallel to the linear hyperplane $v^{\perp} \subseteq F_b$ such that $H$ either halves all measures $\mu_1[b], \dots, \mu_{n}[b]$ or satisfies 
\[
    \mu_i[b](H^+) \le \mu_i[b](F_b)/2 \hspace{3mm} \textrm{ and } \hspace{3mm} \mu_j[b](H^-) \le \mu_j[b](F_b)/2, \hspace{5mm} \text{for some } i,j = 1, \dots,  d+k-1,
\]
where $H^\pm$ is the positive/negative halfspace bounded by $H$ in the direction of $v$.
In either case, if there are multiple such $H$, we choose the middle one.
We write $H = v^{\perp} + a_vv$, for some $a_v \in \R$. Let us denote by 
\begin{equation*}
    w \colon~ \R \longrightarrow \R^{k-1},~ t \longmapsto (t^2, t^3, \dots, t^k)
\end{equation*}
the truncated moment curve and by 
\[
    \gamma^v \colon F_b \longrightarrow \R^{k-1}, \hspace{3mm} x \longmapsto w\left( \langle x, v \rangle -a_v \right)
\]
the \textit{lifting map}, where $\langle x, v \rangle -a_v$ is the signed distance from $x$ to $H$ in the direction of $v$.
Then, the graph of $\gamma^v$ is
\[
    \Gamma^v \coloneqq \{(x, t^2, \dots ,t^k) \in F_b \times \R^{k-1} \colon~ x \in F_b,~ t = \langle x, v \rangle -a_v \}
\]
and we denote by $\pi^v \colon \Gamma^v \cong F_b$ be the natural projection. Then, for each $j=1, \dots, n$, we define the \textit{lift} $\lambda_j^v[b]$ of the measure $\mu_j[b]$ to $\Gamma^v$ to be the pushforward $((\pi^v)^{-1})_*(\mu_j[b])$. See Figure \ref{fig:lifts} for illustration.

The unit vector $n \in S(\R^k)$ has coordinates $(n_1, \dots, n_k)$ and let us denote by 
\[
    n^v \coloneqq (n_1v, n_2, \dots , n_k) \in F_b \times \R^{k-1}.
\]
We define $N^v \subseteq F_b \times \R^{k-1}$ to be an oriented affine hyperplane parallel to $(n^v)^{\perp} \subseteq F_b \times \R^{k-1}$ with a positive halfspace in the direction on $n^v$ which bisects the measure $\lambda_{n}[b]$, that is
\[
    \lambda^v_{n}[b](N^v_+) = \lambda^v_{n}[b](N^v_-).
\]
Again, in case of multiple such hyperplanes $N^v$, we choose the middle one.
Finally, to define the test map \eqref{eq:test-map}, we set $ $
\begin{equation*}
    f_b(v,n) \coloneqq \left(\lambda^v_j[b](N^v_+)-\lambda^v_j[b](N^v_-)\right)_{j=1}^{n-1}.
\end{equation*}
Since all the choices so far have been continuous and the measures $\mu_1[b], \dots, \mu_{n}[b]$ in $F_b$ are absolutely continuous, it follows that $f_b$ is continuous. Moreover, such mapping is continuous in $b \in B$ as well, thus inducing a bundle map $f$ from \eqref{eq:bundle-test-map}.

We now claim that $f_b$ is $\Z_2^2$-equivariant. From the definition it follows that $f_b(v,-n) = - f_b(v,n)$, hence $f$ is equivariant with respect to the second $\Z_2$-factor. On the other hand, if we denote by 
\[
    m \coloneqq (-n_1, n_2, \dots, (-1)^k n_k) \in S\left((\R_{1,2} \oplus \R_2)^{\lfloor k/2 \rfloor} \oplus \R_{1,2}^{k~(\text{mod}~2)}\right)
\]
the vector obtained from $n$ by the action of the first $\Z_2$-factor, we need to show
\begin{equation} \label{eq: f is equivariant, I}
    f_b(v,n) = f_b(-v,m).
\end{equation}
First, we observe that $H = (-v)^{\perp}+(-a_v)(-v)$ stays the same and $a_{-v} = - a_v$. Therefore, the lift map becomes
\begin{equation*}
    \gamma^{-v} \colon F_b \longrightarrow \R^{k-1}, \hspace{3mm} x \longmapsto w\left(\langle x, -v\rangle - a_{-v})\right) = w(-t) = (t^2,-t^3, \dots, (-t)^k),
\end{equation*}
where $t=\langle x, v\rangle - a_{v}$, while its graph is
\begin{equation*}
	\Gamma^{-v} = \{x, t^2, -t^3, \dots , (-t)^k) \in F_b \times \R^{k-1} \colon~ x \in F_b,~ t = \langle x, v \rangle -a_v\}.
\end{equation*}
We define a vector
\[
    m^{-v} \coloneqq (-m_1v, m_2, \dots, m_k) = (n_1v,n_2, -n_3, \dots, (-1)^kn_k) \in F^b \times \R^{k-1}.
\]
Next, we observe that if we write $N^v = \{z\! \in \! F_b \times \R^{k-1} \colon~ \langle n^v,z \rangle = \alpha\}$, for $\alpha \in \R$, then the oriented hyperplane
\begin{equation*}
    M^{-v} \coloneqq \{z \in F_b \times \R^{k-1} \colon~ \langle m^{-v},z \rangle = \alpha\}
\end{equation*}
halves $\lambda_{n}^{-v}[b]$:
\begin{equation} \label{eq:lambda^-v}
	\lambda^{-v}_{n}[b](M^{-v}_+) = \lambda^{-v}_{n}[b](M^{-v}_-).
\end{equation}
Indeed, to see that \eqref{eq:lambda^-v} holds, it is first readily checked that $\langle \gamma^{-v}(y), m^{-v} \rangle = \langle \gamma^{v}(y), n^v \rangle$, which then implies
\begin{align*}
    \pi^{-v}(M^{-v}_+ \cap \Gamma^{-v}) = \{y \in \R^d \colon~ \langle \gamma^{-v}(y), m^{-v} \rangle \ge \alpha \} = \{y \in \R^d \colon~ \langle \gamma^{v}(y), n^v \rangle \ge \alpha\} =  \pi^{v}(N^{v}_+ \cap \Gamma^{v}).
\end{align*}
Therefore, \eqref{eq:lambda^-v} follows from
\begin{align*}\label{eq:lifts+-}
	\lambda^{-v}_{j}[b](M^{-v}_+) &= \mu_{j}[b](\pi^{-v}(M^{-v}_+ \cap \Gamma^{-v}))= \mu_{j}[b](\pi^{v}(N^{v}_+ \cap \Gamma^{v})) = \lambda^{v}_{j}[b](N^{v}_+)
\end{align*}
and from the analogous statement $\lambda^{-v}_{j}(M^{-v}_-) = \lambda^{v}_{j}(N^{v}_-)$, for every $j=1, \dots, n$. Finally, we obtain equivariance \eqref{eq: f is equivariant, I} of the map $f_b$ since
\begin{equation*}
	f_b(-v, m) = \left(\lambda^{-v}_j[b](M^{-v}_+)-\lambda^{-v}_j[b](M^{-v}_-)\right)_{j=1}^{n-1} = \left(\lambda^v_j[b](N^v_+)-\lambda^v_j[b](N^v_-)\right)_{j=1}^{n-1} = f(v,n).
\end{equation*}

\subsection*{Zeros of the test map}

\begin{figure}[h]
        \centering
        \includegraphics[width=.4\linewidth]{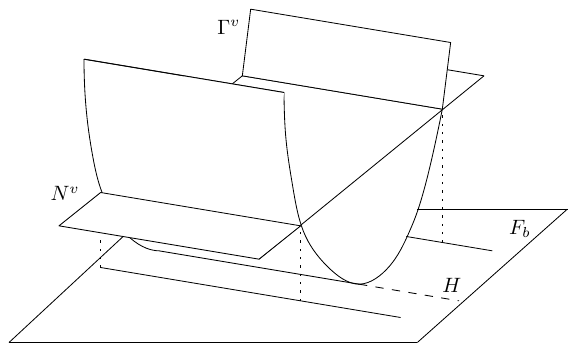}
        \caption{Projection of $\Gamma^v \cap N^v$ is a union of hyperplanes in $F_b$ parallel to $H$.}
        \label{fig:GcapN}
    \end{figure}

Suppose that the test map \eqref{eq:test-map} satisfies $f(v,n)=0$, for some $(v,n) \in S(F_b) \times S(\R^k)$. Then, we have
    \[
        \mu_j[b](\pi^v(N^v_+ \cap \Gamma^v)) = \lambda^v_j[b](N^v_+) = \lambda^v_j[b](N^v_-) = \mu_j[b](\pi^v(N^v_- \cap \Gamma^v)), \hspace{5mm} \text{for }~ j=1, \dots, n,
    \]
    and where $N = \{z \in F_b \times \R^{k-1} \colon~ \langle n^v,z \rangle = \alpha\}$ is an oriented hyperplane in the direction of $n^v$. Let  $p(t) = (a_vn_1-\alpha) + tn_1 + \dots + t^kn_k$ denote a real polynomial in one variable and recall that $H = v^\perp + a_v v$. Then, since $p$ has at most $k$ real zeros, we have that
    \[
        \pi^v(N^v \cap \Gamma^v) = \{x \in F_b \colon~ p(t) = 0,~ t=\langle x,v \rangle -a_v\} = H+\{tv \colon~ t\in \R,~ p(t) = 0\}
    \]
    is a set of at most $k$ affine hyperplanes parallel to $H$, as depicted in Figure \ref{fig:GcapN}. Therefore, the chessboard coloring of $F_b$ with respect to parallel hyperplanes $\pi^v(N^v \cap \Gamma^v)$, given by
    \[
        \pi^v(N^v_{\pm} \cap \Gamma_v) = H+\{tv \colon~ t\in \R,~ \pm p(t) \ge 0\},
    \]
    bisects the masses $\mu_1[b], \dots, \mu_n[b]$.
    
\section{Equivariant maps between bundles}
\label{sec:equiv-bundles}

The main result of the section is Theorem \ref{thm:mass-assign-general} which provides an algebraic condition for the existence of a chessboard bisection of mass assignments on Euclidean vector bundles. As before, notation from Section \ref{sec:prelim} is assumed. Apart from that, for a real vector bundle $\xi \colon \R^d \to E \to B$, we will denote by $w_i(\xi) \in H^i(B;\Z_2)$ the $i$-th Stiefel-Whitney class of $\xi$, where $i=0, \dots, d$. First, we have the following.

\begin{theorem} \label{thm:FH-map-of-bundles}
    Let $d, n \ge 1$ and $a,b \ge 0$ be integers and $\xi \colon~ \R^d \to E \to B$ a Euclidean vector bundle over a paracompact base space with a fiberwise action of $\Z_2^2$ given by the antipodal action of the first factor $\Z_2$. Assume that
    \begin{equation*}
        t_2^{n-1} \notin \left((t_1+t_2)^{a} t_2^b, w_d(\xi) + w_{d-1}(\xi)t_1+\dots+t_1^d\right) \subseteq H^*(B;\FF_2) \otimes \FF_2[t_1, t_2],
    \end{equation*}
    where $|t_1| = |t_2| = 1$.
    Then, there does not exists a $\Z_2^2$-equivariant map of bundles
    \begin{equation} \label{eq:bundle-map}
        \begin{tikzcd}
            S(E) \times S(R^{a,b}) \arrow[r] \arrow[d, "\pi"] &  S(\R_2^{n-1}) \times B \arrow[d, "\mathrm{proj}_2"]\\
            B \arrow[r, equal] & B
        \end{tikzcd}
    \end{equation}
    where $R^{a,b}$ is the $\Z_2^2$-representation defined in Definition $\ref{def:representation-R,a,b}$.
\end{theorem}

Combining Theorem \ref{thm:FH-map-of-bundles} with Lemma \ref{lem:cs-tm-parallel} on the configuration space -- test map scheme, we obtain the following general result about bisections of mass assignments on vector bundles by parallel hyperplanes. Colloquially, it states that if the vector bundle is ``twisted enough", then such a bisection exists.

\begin{theorem}[Bisection of mass assignments by parallel hyperplanes, general]
\label{thm:mass-assign-general}
    Let $d \ge 2$, $k \ge 1$ and $n \ge 1$ be integers. Suppose that $\xi \colon~ \R^d \to E \to B$ is a Euclidean vector bundle over a paracompact base space with a fiberwise action of $\Z_2^2$ given by the antipodal action of the first factor $\Z_2$ and let $\mu_1, \dots, \mu_{n}$ be mass assignments on $\xi$. Assume that
    \begin{equation*}
        t_2^{n-1} \notin \left((t_1+t_2)^{a} t_2^b, w_d(\xi) + w_{d-1}(\xi)t_1+\dots+t_1^d\right) \subseteq H^*(B;\FF_2) \otimes \FF_2[t_1, t_2],
    \end{equation*}
    where $|t_1| = |t_2| = 1$. 
    Then, there exists $b \in B$ such that the masses $\mu_1[b], \dots, \mu_n[b]$ in the fiber $F_b \cong \R^d$ of $\xi$ over $b$ can be simultaneously bisected by a chessboard coloring induced by $k$ or fewer parallel hyperplanes in $F_b$.
\end{theorem}

We now continue with two lemmas needed for the proof of Theorem \ref{thm:FH-map-of-bundles}.

\begin{lemma} \label{lem:proj-bundle-cohomology}
    Let $\xi \colon~ \R^d \to E \to B$ a Euclidean vector bundle over a paracompact base space with a fiberwise antipodal $\Z_2$-action and let $\pi \colon~ S(E) \to B$ denote the projection map of the associated sphere bundle. Then, the map 
    \begin{equation*}
        \pi^* \colon~ H^*_{\Z_2}(B;\FF_2) \longrightarrow H^*_{\Z_2}(S(E);\FF_2)
    \end{equation*}
    is the canonical projection
    \begin{equation*}
        \pi^* \colon~ H^*(B;\FF_2) \otimes \FF_2[t] \longrightarrow H^*(B;\FF_2) \otimes \FF_2[t]/\left(t^d+ t^{d-1}w_1(\xi) + \dots + w_d(\xi)\right),
    \end{equation*}
    where $|t|=1$.
\end{lemma}
\begin{proof}
    Cohomology of the total space of the projectivization $\PP(\xi)$ is represented by
    \begin{equation} \label{eq:cohom-P(E)}
        H^*(\PP(E);\FF_2) \cong H^*(B;\FF_2) \otimes \FF_2[t]/\left(t^d + t^{d-1}w_{1}(\xi) + \dots + w_d(\xi)\right),
    \end{equation}
    where $|t|=1$ is the first Stiefel-Whitney class of the line bundle associated to the twofold cover $S(E) \to \PP(E)$ and the projection $\PP(E) \to B$ induces the inclusion in cohomology (see \cite{Husemoller94fiber}*{Thm.~17.2.5}). Since $\Z_2$ acts freely on $S(E)$, the projection induces a homotopy equivalence
    \[
        S(E) \times_{\Z_2} E\Z_2 \xrightarrow{~\simeq~} S(E)/\Z_2=\PP(E).
    \]
    Therefore, $H^*_{\Z_2}(S(E);\FF_2)$ is isomorphic to \eqref{eq:cohom-P(E)} and the restriction of $\pi^*$ to $H^*(B;\FF_2)$ is the canonical inclusion.
    It is left to show that the class $t$ in \eqref{eq:cohom-P(E)} is the pullback of the generator of $H^*(B\Z_2;\FF_2) \cong \FF_2[t]$ along the projection
    \(
    	\pi_1 \colon~ S(E) \times_{\Z_2} E\Z_2 \to B\Z_2
    \)
    of the Borel fibration. Namely, the two projections from $S(E) \times E\Z_2$ induce a diagram
    \begin{equation*}
        \begin{tikzcd}
            S(E) \arrow[d] & S(E) \times E\Z_2 \arrow[d] \arrow[l, "\simeq"'] \arrow[r] & E\Z_2 \arrow[d] \\
            \PP(E) & \arrow[l, "\simeq"'] S(E) \times_{\Z_2} E\Z_2 \arrow[r, "\pi_1"] & B\Z_2,
        \end{tikzcd}
    \end{equation*}
    where the vertical maps are twofold covers and the two arrows pointing to the left are homotopy equivalences. The claim now follows by naturality of the Stiefel-Whitney classes associated to the line bundles induced by the twofold covers.
\end{proof}

\begin{lemma} \label{lem:FH-equiv-bundle}
    Let $d \ge 1$, $a,b \ge 0$ be integers and $\xi \colon~ \R^d \to E \to B$ a Euclidean vector bundle with a fiberwise action of $\Z_2^2$ given by the antipodal action of the first factor $\Z_2$. Then, for the $\Z_2^2$-representation $R^{a,b}$ defined in Definition $\ref{def:representation-R,a,b}$, the Fadell-Husseini index of the fibration
    \begin{equation*}
        S(\xi) \times S(R^{a,b}) \colon~ S(\R_1^d) \times S(R^{a,b}) \longrightarrow S(E) \times S(R^{a,b}) \xrightarrow{~\pi~} B
    \end{equation*}
    with the fiberwise $\Z_2^2$-action equals
    \begin{align*}
        \FHindex{\Z_2^2}{\FF_2}{\pi} &= \left(t_1^d+t_1^{d-1}w_1(\xi) + \dots +w_d(\xi), (t_1+t_2)^{a} t_2^{b}\right) \subseteq  H^*(B;\FF_2) \otimes \FF_2[t_1, t_2].
    \end{align*}
\end{lemma}
\begin{proof}
    For $i=1,2$, we will denote by $\Z_2^{(i)}$ the $i$'th factor of $\Z_2$ in $\Z_2^2$ and the group cohomology by 
    \[
        H^*_{\Z_2^{(i)}}(\pt;\FF_2) =H^*(B\Z_2^{(i)};\FF_2) \cong \FF_2[t_i], \hspace{3mm} \text{where}~|t_i| = 1.
    \]
    To ease the notation, we set $A \coloneqq H^*(B;\FF_2)$.
    The map $\pi$ in question factors through $\Z_2^2$-equivariant projections as
    \begin{equation*}
        \pi \colon~ S(E) \times S(R^{a,b}) \xrightarrow{~\pi_2~} S(E) \xrightarrow{~\pi_1~} B.
    \end{equation*}
    By Lemma \ref{lem:proj-bundle-cohomology} applied to $\pi_1$ and the action of $\Z_2^{(1)}$, we obtain that 
    \[
        \pi_1^* \colon~ H^*_{\Z_2^{(1)}}(B;\FF_2) \to H^*_{\Z_2^{(1)}}(S(E);\FF_2)
    \]
    is the canonical projection
    \begin{equation} \label{eq:Z2-equiv-cohom-S(E)}
        \pi_1^* \colon~ A[t_1] \longrightarrow A[t_1]/\left(t_1^d+ t_1^{d-1}w_1(\xi) + \dots + w_d(\xi)\right).
    \end{equation}
    In general, for a space $X$ with a $\Z_2^2$-action such that the second factor $\Z_2^{(2)}$ acts trivially, we have 
    \begin{equation*} 
        X \times_{\Z_2^2} E\Z_2^2 = (X \times_{\Z_2^{(1)}} E\Z_2^{(1)}) \times B\Z_2^{(2)}
    \end{equation*}
    and hence
    \begin{equation} \label{eq:equiv-coh-triv-action}
        H^*_{\Z_2^2}(X;\FF_2) \cong H^*_{\Z_2^{(1)}}(X;\FF_2) \otimes \FF_2[t_2]. 
    \end{equation}
    Therefore, since $\Z_2^{(2)}$ acts trivially on the domain and the codomain of $\pi_1$, from \eqref{eq:Z2-equiv-cohom-S(E)} we obtain that 
    \[
        \pi_1^* \colon~ H^*_{\Z_2^{2}}(B;\FF_2) \longrightarrow H^*_{\Z_2^{2}}(S(E);\FF_2)
    \]
    is the canonical projection
    \begin{equation*}
        \pi_1^* \colon~ A[t_1,t_2] \longrightarrow A[t_1,t_2]/\left(t_1^d+ t_1^{d-1}w_1(\xi) + \dots + w_d(\xi)\right).
    \end{equation*}
    We conclude the proof by showing that
    \begin{equation} \label{eq:pi2-in-H-Z2,2}
        \pi_2^* \colon~ H_{\Z_2^2}(S(E); \FF_2) \longrightarrow H^*_{\Z_2^2}(S(E) \times S(R^{a,b});\FF_2)
    \end{equation}
    is the canonical projection
    \begin{align*}
        \pi_2^* \colon~ & A[t_1,t_2]/\left(t_1^d+ t_1^{d-1}w_1(\xi) + \dots + w_d(\xi)\right)\\
        \longrightarrow~ & A[t_1,t_2]/\left(t_1^d+ t_1^{d-1}w_1(\xi) + \dots + w_d(\xi), (t_1+t_2)^{a} t_2^{b}\right).
    \end{align*}
    To do this, let 
    \begin{equation*}
        \eta \colon~ R^{a,b} \longrightarrow (S(E) \times R^{a,b}) \times_{\Z_2^{(1)}} E\Z_2^{(1)} \longrightarrow S(E) \times_{\Z_2^{(1)}} E\Z_2^{(1)}
    \end{equation*}
    be a vector bundle with a fiberwise antipodal $\Z_2^{(2)}$-action and let $\pi_3$ denote the projection of the associated sphere bundle,
    \begin{equation*}
        S(\eta) \colon~ S(R^{a,b}) \longrightarrow (S(E) \times S(R^{a,b})) \times_{\Z_2^{(1)}} E\Z_2^{(1)} \xrightarrow{~\pi_3~} S(E) \times_{\Z_2^{(1)}} E\Z_2^{(1)}.
    \end{equation*}
    The map $\pi_2^*$ from \eqref{eq:pi2-in-H-Z2,2} in $\Z_2^2$-equivariant cohomology is equal to the map induced by $\pi_3$ in $\Z_2^{(2)}$-equivariant cohomology,
    \begin{equation*}
        \pi_2^*=\pi_3^* \colon~ H^*_{\Z_2^{(2)}}(S(E) \times_{\Z_2^{(1)}} E\Z_2^{(1)};\FF_2) \longrightarrow H^*_{\Z_2^{(2)}}((S(E) \times S(R^{a,b})) \times_{\Z_2^{(1)}} E\Z_2^{(1)};\FF_2).
    \end{equation*}
    However, by Lemma \ref{lem:proj-bundle-cohomology} applied to the bundle $\eta$ and by \eqref{eq:Z2-equiv-cohom-S(E)}, we have that
    \begin{align*}
        \pi_3^* \colon~  & A[t_1,t_2]/\left(t_1^d+ t_1^{d-1}w_1(\xi) + \dots + w_d(\xi)\right)\\
        \longrightarrow~ & A[t_1,t_2]/\left(t_1^d+ t_1^{d-1}w_1(\xi) + \dots + w_d(\xi), t_2^{a+b}+t_2^{a+b}w_1(\eta) + \dots + w_{a+b}(\eta)\right)
    \end{align*}
    is the canonical projection. As the final step, we compute
    \begin{equation} \label{eq:sum-of-wi}
    	t_2^{a+b}+t_2^{a+b}w_1(\eta) + \dots + w_{a+b}(\eta) = (t_1+t_2)^{a}t_2^b \in A[t_1,t_2].
    \end{equation}
    Namely, the representation bundle $\rho \colon \R_1 \to \R_1 \times_{\Z_2^{(1)}} E\Z_2^{(1)} \to B\Z_2^{(1)}$ is the canonical line bundle and has the total Stiefel-Whitney class $w(\rho)=1+t_1 \in \FF_2[t_1]$. If we denote by 
    \[
    	\pi_4 \colon~ S(E) \times_{\Z_2^{(1)}} E\Z_2^{(1)} \to B\Z_2^{(1)}
    \]
    the projection of the Borel fibration, we have that the pullback line bundle $\ell \coloneqq \pi_4^*\rho$ may be viewed as
    \begin{equation*}
        \ell \colon~ \R_1 \longrightarrow (S(E) \times \R_1) \times_{\Z_2^{(1)}} E\Z_2^{(1)} \longrightarrow S(E) \times_{\Z_2^{(1)}} E\Z_2^{(1)}.
    \end{equation*}
    Moreover, by Lemma \ref{lem:proj-bundle-cohomology} and naturality of characteristic classes, its total Stiefel-Whitney class is
    \begin{equation} \label{eq:w(ell)}
        w(\ell) = \pi_4^*w(\rho) =  1+t_1 \in A[t_1]/(t_1^d + t_1^{d-1}w_1(\xi) + \dots + w_d(\xi)).
    \end{equation}
    Writing $R^{a,b} = \R_{1,2}^{\oplus a} \oplus \R_2^{\oplus b}$ (see Definition \ref{def:representation-R,a,b}), we have that
    \begin{equation*}
    	\eta \colon~ \R_{1,2}^{\oplus a} \oplus \R_2^{\oplus b} \longrightarrow \left(S(E) \times (\R_{1,2}^{\oplus a} \oplus \R_2^{\oplus b})\right) \times_{\Z_2^{(1)}} E\Z_2^{(1)} \longrightarrow S(E) \times_{\Z_2^{(1)}} E\Z_2^{(1)}
    \end{equation*}
    is isomorphic to the direct sum of $\ell^{\oplus a}$ and the rank $b$ trivial bundle. Therefore, by \eqref{eq:w(ell)} and the Whitney sum formula, we have
    \begin{equation*}
    	w(\eta) = w(\ell^{\oplus a}) = (1+t_1)^a \in A[t_1]/(t_1^d + t_1^{d-1}w_1(\xi) + \dots + w_d(\xi)).
    \end{equation*}
	Finally, this implies \eqref{eq:sum-of-wi}, since now
	\begin{align*}
		t_2^{a+b}+t_2^{a+b}w_1(\eta) + \dots + w_{a+b}(\eta) = (t_2^a + t_2^{a-1}w_{1}(\eta) + \dots + w_a(\eta)) t_2^b = (t_1+t_2)^at_2^b,
	\end{align*}
	which completes the proof.
\end{proof}

Finally, we may prove the main result of the section.

\begin{proof}[Theorem \ref{thm:FH-map-of-bundles}]
    We will prove the contrapositive. Assuming that the $\Z_2^2$-equivariant bundle map \eqref{eq:bundle-map} exists, by monotonicity of the Fadell-Husseini index, we have that
    \[
        \FHindex{\FF_2}{\Z_2^2}{\mathrm{proj}_2} \subseteq \FHindex{\FF_2}{\Z_2^2}{\pi}.
    \]
    The index of $\pi$ is computed in Lemma \ref{lem:FH-equiv-bundle}. On the other hand, applying Lemma \ref{lem:proj-bundle-cohomology} to the sphere bundle $\mathrm{proj}_2$ coming from the trivial vector bundle, we get that
    \begin{equation*}
        \mathrm{proj}_2^* \colon~H^*_{\Z_2^{(2)}}(B;\FF_2) \longrightarrow H^*_{\Z_2^{(2)}}(S(\R_2^{n-1}) \times B;\FF_2)
    \end{equation*}
    is the canonical projection
    \begin{equation*}
        \mathrm{proj}_2^* \colon~H^*(B;\FF_2) \otimes \FF_2[t_2] \longrightarrow H^*(B;\FF_2) \otimes \FF_2[t_2]/(t_2^{n-1})
    \end{equation*}
    in the equivariant cohomology with respect to the second factor $\Z_2^{(2)}$ in $\Z_2^2$.
    Similarly as in the proof of Lemma \ref{lem:FH-equiv-bundle} (see  \eqref{eq:equiv-coh-triv-action}), since the first factor $\Z_2^{(1)}$ acts trivially on the total space of $\mathrm{proj}_2$, we conclude that the map in $\Z_2^2$-equivariant cohomology
    \begin{equation*}
        \mathrm{proj}_2^* \colon~H^*_{\Z_2^{2}}(B;\FF_2) \longrightarrow H^*_{\Z_2^{2}}(S(\R_2^{n-1}) \times B;\FF_2)
    \end{equation*}
    is the canonical projection
    \begin{equation*}
        \mathrm{proj}_2^* \colon~H^*(B;\FF_2) \otimes \FF_2[t_1,t_2] \longrightarrow H^*(B;\FF_2) \otimes \FF_2[t_1,t_2]/(t_2^{n-1}).
    \end{equation*}
    Therefore, $\FHindex{\FF_2}{\Z_2^2}{\mathrm{proj}_2} = (t_2^{n-1})$, which finishes the proof of the theorem.
\end{proof}


\section{Proof of the main result}
\label{sec:proof-main}

In this section, we prove the main result of the paper. First, we recall its statement.

\thmMain*

We have the following technical result needed for the proof of the theorem. We use the notation for the cohomology of the Grassmannian described in \eqref{eq:coh-grass} in Section \ref{sec:prelim}.

\begin{lemma} \label{lem:ideal-stirling}
	Let $d \ge 1$, $k \ge 1$ and $m \ge 0$ be integers. Then,
	\begin{equation} \label{eq:ideal-non-membership}
			t_2^{d+m+k-2} \notin \left(t_1^d+t_1^{d-1}w_1 + \dots +w_d, ~ (t_1+t_2)^{\lceil k/2 \rceil} t_2^{\lfloor k/2 \rfloor}\right) \subseteq H^*(\grassR{d}{d+m} ;\FF_2) \otimes \FF_2[t_1, t_2]
	\end{equation}
	if and only if the Stirling number of second kind $S(d+m+k-1,k)$ is odd.
\end{lemma}

With the lemma at hand, whose proof is deterred to the end of the section, we may now prove our main result.

\begin{proof} [Proof of Theorem \ref{thm:assignment}]
For $n \coloneqq d+m+k-1$, let $\mu_1, \dots, \mu_n$ be the mass assignments on the tautological vector bundle $\R^d \to \gamma_d \to \grassR{d}{d+m}$. By Theorem \ref{thm:mass-assign-general}, there exists a linear $d$-plane $V \in \grassR{d}{d+m}$ such that the measures $\mu_1[V], \dots, \mu_n[V]$ in $V$ can be bisected by the chessboard coloring of at most $k$ parallel hyperplanes in $V$ if
\begin{equation*}
			t_2^{n-1} \notin \left(t_1^d+t_1^{d-1}w_1 + \dots +w_d, ~ (t_1+t_2)^{\lceil k/2 \rceil} t_2^{\lfloor k/2 \rfloor}\right) \subseteq H^*(\grassR{d}{d+m} ;\FF_2) \otimes \FF_2[t_1, t_2].
\end{equation*}
By Lemma \ref{lem:ideal-stirling}, the latter is equivalent to $S(d+m+k-1,k)$ being odd.
\end{proof}

As the final step, we provide the deterred proof of the lemma.

\begin{proof} [Proof of Lemma \ref{lem:ideal-stirling}]
	Let $A \coloneqq H^*(\grassR{d}{d+m} ;\FF_2)$. Then, the ambient ring for \eqref{eq:ideal-non-membership} is the free $A$-algebra $A[t_1,t_2]$, where $|t_1| = |t_2| = 1$.
Contrapositive of \eqref{eq:ideal-non-membership} is equivalent to:
\begin{equation} \label{eq:ideal-1st}
	\left( \exists p,q \in A[t_1, t_2]\right) \hspace{3mm} t_2^{d+m+k-2} = p \cdot \Big( \sum_{i=0}^d t_1^{d-i}w_i \Big) +  q \cdot (t_1+t_2)^{\lceil k/2 \rceil} t_2^{\lfloor k/2 \rfloor} \in A[t_1,t_2].
\end{equation}
The proof continues by showing, in several steps, that \eqref{eq:ideal-1st} is equivalent to $S(d+m+k-1,k)$ being even. Before proceeding further, we note that the relations \eqref{eq:cohom-grass} in $A$ imply the identity
\begin{align}\label{eq:product-identity}
	( t_1^d+t_1^{d-1}w_1 + \dots + w_d ) \left( t_1^m + t_1^{m-1}\overline{w}_1 + \dots + \overline{w}_m \right) = t_1^{d+m} \in A[t_1,t_2],
\end{align}
which will be used in several occasions in the proof.

\noindent
\textbf{Step 1:} \eqref{eq:ideal-1st} is equivalent to
\begin{equation} \label{eq:ideal-2nd}
	\left( \exists p,q \in A[t_1, t_2]\right) \hspace{3mm} t_2^{d+m+\lceil k/2 \rceil-2} = p \cdot \Big( \sum_{i=0}^d t_1^{d-i}w_i \Big) +  q \cdot (t_1+t_2)^{\lceil k/2 \rceil} \in A[t_1,t_2].
\end{equation}
Indeed, \eqref{eq:ideal-1st} is obtained by multiplying both sides of \eqref{eq:ideal-2nd} by $t_2^{\lfloor k/2 \rfloor}$ and using the fact that $k = \lceil k/2 \rceil + \lfloor k/2 \rfloor$. For the other implication, if $p$ in \eqref{eq:ideal-1st} was divisible by $t_2^{\lfloor k/2 \rfloor}$, then we could divide both sides of the equality with this monomial to obtain \eqref{eq:ideal-2nd}. To see that the divisibility assumption on $p$ holds, we first notice that from \eqref{eq:ideal-1st} it follows that $p \cdot ( t_1^d+t_1^{d-1}w_1 + \dots +w_d )$ is divisible by $t_2^{\lfloor k/2 \rfloor}$. Hence, by the identity \eqref{eq:product-identity}, we have that
\begin{equation*}
	p  \left( t_1^d+t_1^{d-1}w_1 + \dots + w_d \right) \left( t_1^m + t_1^{m-1}\overline{w}_1 + \dots + \overline{w}_m \right) = pt_1^{d+m} \in A[t_1, t_2]
\end{equation*}
is divisible by $t_2^{\lfloor k/2 \rfloor}$, and so must be $p$, because $t_1$ and $t_2$ are $A$-algebra generators of $A[t_1,t_2]$.

\noindent
\textbf{Step 2:} \eqref{eq:ideal-2nd} is equivalent to
\begin{equation} \label{eq:ideal-3rd}
	\left( \exists p,q \in A[t_1, t_2]\right) \hspace{3mm} (t_1+t_2)^{d+m+\lceil k/2 \rceil-2} = p \cdot \Big( \sum_{i=0}^d t_1^{d-i}w_i \Big) +  q \cdot t_2^{\lceil k/2 \rceil} \in A[t_1,t_2].
\end{equation}
Indeed, an $A$-algebra involution
\begin{equation*}
	A[t_1,t_2] \longrightarrow A[t_1,t_2],~t_1 \longmapsto t_1,~t_2 \longmapsto t_1+t_2,
\end{equation*} 
transforms \eqref{eq:ideal-2nd} to \eqref{eq:ideal-3rd}, and vice versa.

\noindent
\textbf{Step 3:} \eqref{eq:ideal-3rd} is equivalent to
\begin{equation} \label{eq:ideal-4th}
	\left( \exists p,q \in A[t_1, t_2]\right) \hspace{3mm} t_1^{d+m-1}t_2^{\lceil k/2 \rceil-1}  \cdot S = p \cdot \Big( \sum_{i=0}^d t_1^{d-i}w_i \Big) +  q \cdot t_2^{\lceil k/2 \rceil} \in A[t_1,t_2],
\end{equation}
where $S \coloneqq \binom{d+m+\lceil k/2 \rceil-2}{\lceil k/2 \rceil-1} \in \Z$. Indeed, in the binomial formula
\[
	(t_1+t_2)^{d+m+\lceil k/2 \rceil-2} = \sum_{i+j=d+m+\lceil k/2 \rceil-2} t_1^it_2^j \tbinom{d+m+\lceil k/2 \rceil-2}{i},
\]
all of the summands with $(i,j) \neq (d+m-1, \lceil k/2 \rceil-1)$ have either $i \ge d+m$ or $j \ge \lceil k/2 \rceil$, hence are divisible by $t_1^{d+m}$ or $t_2^{\lceil k/2 \rceil}$. Therefore, by \eqref{eq:product-identity}, they are in the ideal
\begin{equation*}
	 \left( t_1^{d+m}, t_2^{\lceil k/2 \rceil}\right) \subseteq \left(t_1^d+t_1^{d-1}w_1 + \dots +w_d,~ t_2^{\lceil k/2 \rceil} \right) \subseteq A[t_1, t_2].
\end{equation*}

\noindent
\textbf{Step 4:} \eqref{eq:ideal-4th} is equivalent to $S(d+m+k-1,k)$ being even. To show this, we first note that $S(d+m+k-1, k)$ and $S$ have the same parity (see \cite{Stanley11}*{(1.94c)} or \cite{ChanManna10}*{Thm.~2.1}). If $S$ is even, we have that \eqref{eq:ideal-4th} holds trivially. On the other hand, assume that \eqref{eq:ideal-4th} holds. By multiplying both sides of the equation with $t_1^m + t_1^{m-1}\overline{w}_1 + \dots + \overline{w}_m$, we obtain via the identity \eqref{eq:product-identity}, that
\begin{equation*}
	t_1^{d+m-1}t_2^{\lceil k/2 \rceil-1} \cdot S \cdot \left( t_1^m + t_1^{m-1}\overline{w}_1 + \dots + \overline{w}_m \right) = p\cdot t_1^{d+m} + q' \cdot t_2^{\lceil k/2 \rceil} \subseteq A[t_1, t_2],
\end{equation*}
after adjusting $q$ to $q'$. However, if we exclude the last summand on the left hand side, the rest of the sum is divisible by $t_1^{d+m}$. Thus, after adjusting $p$ to $p'$, we get
\begin{equation*}
	t_1^{d+m-1}t_2^{\lceil k/2 \rceil-1} \cdot S \cdot \overline{w}_m  = p'\cdot t_1^{d+m} + q' \cdot t_2^{\lceil k/2 \rceil} \subseteq A[t_1, t_2],
\end{equation*}
which implies $S \cdot \overline{w}_m = 0 \in A$. Therefore, $S$ is even, since $\overline{w}_m \neq 0 \in A$ (which follows, for example, from the relations \eqref{eq:cohom-grass}). This completes the proof. 
\end{proof}

\bibliographystyle{plain}
\bibliography{references.bib}

\end{document}